\documentclass{article}
\usepackage{amssymb,graphpap,xcolor,xypic,color}
\usepackage{graphicx}
      \usepackage{amsmath}
      \usepackage{amssymb}
      \usepackage{amsfonts}
      \usepackage{epsf}
      \usepackage{color}
      \usepackage{enumerate}

\usepackage{geometry}
\usepackage{doc}

\usepackage{url}

\usepackage{graphicx}
\usepackage{epstopdf}
\title{Introduction to Noncommutative Algebraic Geometry}
\author{Manizheh Nafari}
\date{}

\begin{document}

\maketitle

\abstract{
This Lecture Notes is meant to introduce noncommutative algebraic geometry tools (which were invented by M. Artin, W. Schelter, J. Tate, and M. Van den Bergh in the late 1980s) and also graded skew Clifford algebras (which were introduced by T. Cassidy and M. Vancliff).}

\pagebreak
\tableofcontents



\pagebreak

%
%
\section{Introduction}
\label{introduction}

M. Artin, W. Schelter, J. Tate, and M. Van den Bergh introduced the notion of non-commutative regular algebras and invented new methods in algebraic geometry in the late 1980s to study them (\cite{AS}, \cite{ATV1}, \cite{ATV2}). Such algebras are viewed as non-commutative analogues of polynomial rings; indeed, polynomial rings are examples of regular algebras.\\

By the 1980s, a lot of algebras had arisen in quantum physics, specifically quantum groups, and many traditional algebraic techniques failed on these new algebras. In physics, quantum groups are viewed as algebras of non-commuting functions acting on some ``non-commutative space''(\cite{D}). In the early 1980s, E. K. Sklyanin, a physicist, constructed a family of graded algebras on four generators (\cite{S}). These algebras were later proved to depend on an elliptic curve and an automorphism (\cite{FO}). By the late 1980s, it was known that many of the algebras in quantum physics are regular algebras; in particular, the family of algebras constructed by Sklyanin consists of regular algebras.\\

The main results in \cite{AS}, \cite{ATV1}, and \cite{ATV2} concern the classification of regular algebras of global dimension 3 on degree-one generators. The quadratic regular algebras of global dimension 3 can be described using geometry, i.e. the point scheme $E \subseteq \mathbb{P}^2$. These algebras, where $E$ contains a line as well as those that are ``generic'', are given in \cite{ATV1}, and \cite{ATV2}, and entail: $\mathbb{P}^2$, elliptic curve, conic union a line, triangle, (triple) line, a union of $n$ lines where $n \in \{2,3\}$ with one intersection point. It should be noted that the cases where $E$ is a nodal cubic curve or a cuspidal cubic curve are not discussed in \cite{ATV1} or \cite{ATV2} as such algebras are not generic.\\

T. Cassidy and M. Vancliff introduced a class of algebras that provide an ``easy'' way to write down some quadratic regular algebras of global dimension $n$  where $n \in \mathbb{N}$ (\cite{CV}). In fact, they generalized the notion of a graded Clifford algebra and called it a graded skew Clifford algebra (see Definition 3.1).\\

\section{Definitions}

\subsection{Definition of Graded Algebras \cite{ATV1}}

Throughout this lecture notes, $\mathbb{K}$ denotes an algebraically closed field, $\text{char}(\mathbb{K}) \neq 2$, and $\mathbb{K}^{\times}$ denotes $\mathbb{K} \setminus \{ 0 \}$.

\noindent A $\mathbb{K}$-algebra $A$ is called a graded algebra if:

\begin{enumerate}
\item[(1)] $ A = \bigoplus_{i \geq 0} A_i$ where the $A_i$ are vector spaces over $\mathbb{K}$,
\item[(2)] $ \text{dim} A_1 < \infty $,
\item[(3)] $ A_i A_j \subseteq A_{i + j} $ for all $i,j$,
\item[(4)] $A _0 = \mathbb{K}$,
\item[(5)] $A$ generated by $A_1$ only.
\end{enumerate}

\noindent For each $i$, $A_i$ is the span of the homogeneous elements of degree $i$.

\subsection{Examples}

\noindent (1) The polynomial ring $A = \mathbb{K}[x_1, \ldots, x_d]$ where $x_1,\ldots,x_d$ have degree 1.\\
Here,
\[
A_1 = \mathbb{K} x_1 \oplus \mathbb{K} x_2 \oplus \cdots \oplus \mathbb{K} x_d,
\]
and
\[
\text{dim}_{\mathbb{K}} A_i = \binom{i + d - 1}{d - 1} \quad \text{for all} \quad i \quad \text{(c.f., \cite{MR})}.
\]

\noindent (2) The free algebra $A = \mathbb{K} \langle x_1,\ldots,x_d \rangle$ where $x_i$, for all $i$, have degree $n_i \in \mathbb{Z}$.\\
Here, $A$ is a non-commutative analogue of the algebra $A$ in (1).

\subsection{Nonexamples}

\noindent (1) The algebra
\[
A = \frac{\mathbb{K}[x,y]}{\langle x^2 - y \rangle},
\]
where $x$ and $y$ have degree 1, is not graded. The relation $x^2 = y$ in $A$ is not homogeneous and so $A_1 \cap A_2 \neq \{ 0 \}$ which violates (1) in Definition 2.1.1.\\

\noindent (2) The algebra
\[
A = \frac{\mathbb{K}[x,y]}{\langle x^2 - y \rangle},
\]
where $x$ has degree 1 and $y$ has degree 2, is graded but not generated by $A_1$ since $y \in A_2$.

\subsection{Definition of Quadratic $\mathbb{K}$-Algebra}

\noindent A $\mathbb{K}$-algebra $A$ is called quadratic if:

\begin{enumerate}
\item[(1)] $A$ is graded (as defined above),
\item[(2)] $A$ is a quotient of the free algebra by homogeneous relations of degree 2.
\end{enumerate}

\subsection{Example}

\noindent The algebra
\[
\mathbb{K}[x_1, \ldots, x_d] = \frac{\mathbb{K} \langle x_1,\ldots,x_d \rangle}{\langle x_i x_j - x_j x_i; 1 \leq i,j \leq d \rangle}, \quad \text{deg}(x_i) = 1 \quad \text{for all} \quad i
\]
is quadratic.

\subsection{Nonexample}

\noindent The algebra
\[
A = \frac{\mathbb{K}[x]}{\langle x^3 \rangle}, \quad \text{where} \quad x \quad \text{has degree 1},
\]
is graded but is not quadratic. The relation $x^3 = 0$ has degree 3.\\

In order to define a regular algebra, we first need the concepts of polynomial growth, global dimension, and Gorenstein, which we now define.

\subsection{Global Dimension}

\noindent The algebra $A$ has global dimension $d < \infty$ if every $A$-module $M$ has projective dimension $\leq d$ and there exists at least one module $M$ with projective dimension $d$.

\subsection{Example}

\noindent The polynomial ring, $\mathbb{K}[x_1,\ldots, x_d]$, has global dimension $d$ by Hilbert's syzygy theorem (c.f., \cite{R}).

\subsection{Definition of Polynomial Growth (c.f.,\cite{MR})}

\noindent A graded algebra $A$, as above, is said to have polynomial growth if there exists positive real numbers $ c, \delta$ such that
\[
\text{dim}_{\mathbb{K}} A_n \leq c n^\delta \quad \text{for all} \quad n \gg 0.
\]
For all known quadratic regular algebras of global dimension $d$, the minimal such $\delta$ is $d-1$ (\cite[$\S$2]{ATV1}).

\subsection{Example}

\noindent Let $A = \mathbb{K}[x_1, x_2]$, then
\[
\text{dim}_{\mathbb{K}} A_n = \binom{n + 1}{1} = n + 1 \leq n^{1 + \epsilon},
\]
for all $\epsilon > 0$ where $n \gg 0$. Thus $A$ has polynomial growth.

\subsection{Definition of Gorenstein \cite{AS}}

\noindent By \cite[$\S$2]{ATV1}, for a graded algebra $A$ as in Definition 2.1.1, the global dimension of $A$ equals the projective dimension of the graded left module $_ A \mathbb{K}$ (and projective dimension of the right module $\mathbb{K}_A$).\\
\noindent The algebra $A$ is Gorenstein if
\begin{enumerate}
\item[(1)] the projective modules $P^i$ appearing in a minimal resolution
\[
0 \rightarrow P^d \rightarrow ... \rightarrow P^1 \rightarrow P^0 \rightarrow _ A \mathbb{K} \rightarrow 0
\]
of $_A \mathbb{K}$ are finitely generated, and if
\item[(2)] applying the functor
\[
M \rightsquigarrow M^{*} : = Hom_{A} (M,A) = \{ \text{graded homomorphisms}: M \rightarrow A \}
\]
to the resolution in (1) yields a projective resolution
\[
0 \rightarrow P^{0*} \rightarrow P^{1*} \rightarrow ... \rightarrow P^{d*} \rightarrow \mathbb{K}_{A} \rightarrow 0
\]
of the graded right $A$-module $\mathbb{K}_A$.
\end{enumerate}

\subsection{Example}

\noindent The algebra
\[
A = \frac{\mathbb{K} \langle x,y \rangle}{\langle x y - q y x \rangle}, \quad \text{where} \quad q \in \mathbb{K}^{\times},
\]
is Gorenstein (\cite[$\S$0]{AS}).

\subsection{Definition of Regular Algebras \cite{ATV1}}

\noindent A graded $\mathbb{K}$-algebra $A$ is called a regular algebra if

\begin{enumerate}
\item[(1)] $A$ has polynomial growth,
\item[(2)] $A$ has finite global dimension,
\item[(3)] $A$ is Gorenstein.
\end{enumerate}

\subsection{Definition of Normalizing Sequence}

\noindent A sequence $a_1,\ldots,a_n$ of elements of a ring $R$ with identity is called a normalizing sequence if $a_1$ is normal element in $R$ (i.e. $a_1 R = R a_1$) and for each $j \in \{ 1,\ldots, n-1 \}$, $a_{j+1}$ is a normal element in $R / \sum_{i=1}^{j} a_{i}R$ and also $\sum_{i=1}^{n} a_{i}R \neq R.$

\section{Graded Skew Clifford Algebras}

T. Cassidy and M. Vancliff defined a class of algebras in \cite{CV} that provide an ``easy'' way to write down some quadratic regular algebras of global dimension $d$ for all $d \in \mathbb{N}$.

\subsection{Definition of Graded Skew Clifford Algebras \cite{CV}}

\noindent For $\{i, j \}\subset \{ 1, \ldots , n\}$, let $\mu_{ij} \in \mathbb{K}^\times$ satisfy $\mu_{ij} \mu_{ji} = 1$ for all $i \neq j$, and write $\mu = (\mu_{ij}) \in M(n, \mathbb{K})$. A matrix $M \in M(n, \mathbb{K})$ is called $\mu$-symmetric if $M_{ij} = \mu_{ij}M_{ji}$ for all $i, j = 1, \ldots , n$.\\
\noindent Henceforth, suppose $\mu_{ii} = 1$ for all $i$, and fix $\mu$-symmetric matrices $M_1, \ldots , M_n \in M(n, \mathbb{K})$. A graded skew Clifford algebra associated to $\mu$ and $M_1$, $\ldots ,$ $M_n$ is a graded $\mathbb{K}$-algebra on degree-one generators $x_1, \ldots , x_n$ and on degree-two generators $y_1, \ldots , y_n$ with defining relations given by:
\begin{enumerate}
\item[(a)] $ x_i x_j + \mu_{ij} x_j x_i = \sum_{k=1}^n (M_k)_{ij} y_k$ for all $i, j = 1, \ldots , n$, and
\item[(b)] the existence of a normalizing sequence $\{ r_1, \ldots , r_n \}$ of homogeneous elements that span $\mathbb{K} y_1 + \cdots + \mathbb{K} y_n$.
\end{enumerate}

\subsection{Example}

\noindent Let $\mu_{21},\lambda \in \mathbb{K}^{\times}$. If
\[
M_1 = \left[\begin{array}{cc} 0 & 1 \\ \mu_{21} & 0 \end{array}\right], \qquad M_2 = \left[\begin{array}{cc} 2 & 0 \\ 0 & 2 \lambda \end{array}\right],
\]
then any graded skew Clifford algebra $A$ associated to $M_1, M_2$ satisfies
\[
\frac{\mathbb{K} \langle x_1,x_2 \rangle}{\langle {x_2}^2 - \lambda {x_1}^2 \rangle} \twoheadrightarrow A
\]
since
\[
x_1 x_2 + \mu_{12} x_2 x_1 = y_1, \qquad y_2 = {x_1}^2, \qquad \lambda y_2 = {x_2}^2.
\]

\subsection{Definition of Quadric System \cite{CV}}

\noindent Let $S$ be the $\mathbb{K}$-algebra on generators $z_1,\ldots,z_n$ with defining relations
\[
z_j z_i = \mu_{ij} z_i z_j, \qquad \text{for all} \quad i,j
\]
and let
\[
q_k: = \left[\begin{array}{ccc} z_1 & \ldots & z_n\\ \end{array}\right] \quad M_k \quad \left[\begin{array}{ccc} z_1 \\ \vdots \\ z_n \end{array}\right] \quad \in S.
\]
We say $\{ q_1,\dots,q_n \}$ is a quadric system.

\subsection{Example}

\noindent For the algebra $A$ in Example 2.2.2, we have
\[
S = \frac{\mathbb{K} \langle z_1,z_2 \rangle}{\langle z_2 z_1 - \mu_{12} z_1 z_2 \rangle}.
\]
Moreover,
\[
q_1 = 2 z_1 z_2, \quad q_2 = 2 {z_1}^2 + 2 \lambda {z_2}^2.
\]
However, since $\text{char}(\mathbb{K}) \neq 2$, we consider:
\[
q_1 = z_1 z_2, \quad q_2 = {z_1}^2 + \lambda {z_2}^2.
\]

\subsection{Definition of Normalizing Quadric System}

\noindent A quadric system $\{ q_1,\dots,q_n \}$ is normalizing if $\sum_{k=1}^{n} \mathbb{K} q_k \subset S$ is spanned by a normalizing sequence of $S$.

\subsection{Example}

\noindent Referring to Example 2.2.4, in $S$, $z_i$ is normal for all $i$, and
\[
q_1 z_1 = \mu_{12} z_1 q_1, \quad q_1 z_2 = \mu_{21} z_2 q_1 .
\]
Therefore $q_1$ is normal in $S$.\\
In $\frac{S}{\langle q_1 \rangle}$, we have
\[
q_2 z_1 = z_1 ({z_1}^2 + \lambda {\mu_{12}^2} {z_2}^2), \quad q_2 z_2 = {\mu_{21}}^2 z_2 ( {z_1}^2 + \lambda {\mu_{12}}^2 {z_2}^2).
\]
So $q_2$ is normal in $\frac{S}{\langle q_1 \rangle}$ if $\lambda = 0$ or if $\lambda \neq 0$ and ${\mu_{12}}^2 = 1$.

\subsection{Definition of Zero Locus \cite{CV}}

\noindent Suppose $A = \mathbb{K} \langle x_1, \ldots, x_n \rangle$ and $f \in A_2$. We define the zero locus $\mathcal V(f)$ of $f$ to be
\[
\mathcal V(f) = \{ p \in \mathbb{P}^ {n - 1} \times \mathbb{P}^{n - 1}: f(p) = 0 \},
\]
where $\mathbb{P}^{n-1}$ is identified with $\mathbb{P}(A_1^*)$.\\

\noindent Similarly if $f_1,\ldots,f_m \in A_2$, then
\[
\mathcal V(f_1, \ldots, f_m) = \{ p \in \mathbb{P}^ {n - 1} \times \mathbb{P}^{n - 1}: f_i(p) = 0 \quad \text{for all} \quad i \}.
\]

\subsection{Definition of Base-Point Free \cite{CV}}

\noindent Let $Z$ be the zero locus in $\mathbb{P}^{n-1} \times \mathbb{P}^{n-1}$ of the defining relations of $S$, i.e.
\[
Z = \bigcap_{i,j} \mathcal V(z_j z_i - \mu_{ij} z_i z_j) \subset \mathbb{P}^{n-1} \times \mathbb{P}^{n-1}.
\]
The quadric system $\{ q_1,\ldots,q_n \}$ is said to be base-point free (BPF) if $Z \cap \mathcal V(q_1,\ldots,q_n)$ is empty.

\subsection{Example}

\noindent Referring to Example 2.2.4, let
\[
p = ((\alpha_1, \alpha_2),(\beta_1, \beta_2)) \in \mathbb{P}^1 \times \mathbb{P}^1,
\]
and let
\[
(z_2 z_1 - \mu_{12} z_1 z_2) (p) = 0.
\]
Therefore, we have
\[
\alpha_2 \beta_1 - \mu_{12} \alpha_1 \beta_2 = 0.
\]
If $\alpha_2 = 0$, then $\beta_2 = 0$. So $((1,0),(1,0)) \in \mathbb{P}^1 \times \mathbb{P}^1$. \\
If $\alpha_2 \neq 0$, i.e., $\alpha_2 = 1$, then $\beta_1 = \mu_{12} \alpha_1 \beta_2$. So, $((\alpha_1,1),(\mu_{12} \alpha_1,1)) \in \mathbb{P}^1 \times \mathbb{P}^1$. \\
Therefore,
\[
Z = \{ ((\alpha_1,\alpha_2),(\mu_{12} \alpha_1, \alpha_2)): (\alpha_1,\alpha_2) \in \mathbb{P}^1 \}.
\]
Let $p \in Z$. We have
\[
0 = q_1 (p) = \alpha_1 \alpha_2, \quad 0 = q_2(p) = \mu_{12} {\alpha_1}^2 + \lambda {\alpha_2}^2.
\]
Thus $\alpha_1 = \alpha_2 = 0$ which is contradiction. Therefore $\{ q_1,q_2 \}$ is BPF.


\begin{thebibliography}{99}


\bibitem{AS}
Artin, M. and Schelter, W., \emph{Graded Algebras of Global Dimension 3}, Adv. Math., 66 (1987), 171-216.

\bibitem{ATV1}
M. Artin, J. Tate and M. Van den Bergh, \emph{Some Algebras Associated to Automorphisms of Elliptic Curves}, The Grothendieck Festschrift 1, Eds. P.
Cartier et al. Birkhauser (1990), 33-85.

\bibitem{ATV2}
Artin, M., Tate, J., and Van den Bergh, M., \emph{Modules Over Regular Algebras of Dimension 3}, Invent. Math., 106 (1991), 335-388.

\bibitem{CV}
Cassidy, T. and Vancliff, M., \emph{Generalizations of Graded Clifford Algebras and of Complete Intersections}, Journal of the London Mathematical
Society 81 (2010), 91-112.

\bibitem{D}
Drinfel'd, V. G., \emph{Quantum Groups}, Proc. Int. Cong. Math., Berkeley 1 (1986), 798-820.


\bibitem{FO}
Feigin, B. L., and Odesskii, A. B., \emph{Elliptic Sklyanin Algebras}, Func. Anal. Appl. 23 (1989), 45-54.





\bibitem{MR}
McConnell, J. C. and Robson, J. C., \emph{Noncommutative Noetherian Rings}, Graduate Studies in Mathematics, American Mathematical Society, 2001.


\bibitem{R}
Rotman, Joseph J., \emph{An Introduction to Homological Algebra}, Second edition, Universitext, Springer, New York, 2009, xiv+709 pp.

\bibitem{S}
Sklyanin, E. K., \emph{Some Algebraic Structures Connected to the Yang-Baxter Equation}, Func. Anal. Appl. 16 (1982), no. 4, 27-34.




\end{thebibliography}

\end{document}